\documentclass[11pt,a4paper]{amsart}
\usepackage[T1]{fontenc}
\usepackage[utf8]{inputenc}
\usepackage{lmodern}
\usepackage[margin=1in]{geometry}
\usepackage{amsmath,amssymb,amsthm}
\usepackage{booktabs,array}
\usepackage{microtype}
\usepackage[hidelinks]{hyperref}
\usepackage{enumitem}
\usepackage{listings}
\newtheorem{theorem}{Theorem}[section]
\newtheorem{lemma}[theorem]{Lemma}
\newtheorem{proposition}[theorem]{Proposition}
\newtheorem{corollary}[theorem]{Corollary}
\theoremstyle{definition}
\newtheorem{definition}[theorem]{Definition}
\newtheorem{example}[theorem]{Example}
\theoremstyle{remark}
\newtheorem{remark}[theorem]{Remark}
\DeclareMathOperator{\Cay}{Cay}

\newcommand{\Pow}{\mathcal{P}}
\newcommand{\Comp}{\mathcal{C}}
\newcommand{\gen}[1]{\langle #1\rangle}

\title{Connected Cayley Approximation of Power Graphs of Finite Groups}
\date{}
\begin{document}
\maketitle
\begin{center}
\textbf{Zahra Rabiei}\\
{\small Department of Mathematical Sciences, Yazd University, Yazd, Iran,\\ \texttt{123456@stu.ac.ir}}\par
\medskip
\textbf{Mohammad A. Iranmanesh}\textsuperscript{*}\\
{\small Department of Mathematical Sciences, Yazd University, Yazd, Iran}\par
\smallskip
{\small \textsuperscript{*}Corresponding author: \texttt{iranmanesh@yazd.ac.ir}}
\end{center}
\begin{abstract}
For a nontrivial finite group \(G\), we define its connected intersection density as the maximum, over all connected Cayley graphs on \(G\), of the
proportion of edges that also belong to the undirected power graph of \(G\). We show that this density is \(1\) if and only if \(G\) is cyclic, and that
it is \(1/4\), \(1/6\), and \(1/12\) for \(A_4\), \(S_4\), and \(A_5\), respectively. For the dihedral group of order \(2n\), with \(n\geq 3\), we give a general formula and show that the connected intersection density equals \((\binom{n}{2}+1)/n^2\) if and only if \(n\) is a prime power. We also give a formula for the connected intersection density of dicyclic groups. For generalized quaternion groups of order \(2N\), where \(N=2^{m-1}\geq 4\), the density is \(\frac{N^2+8}{2N(N+1)}\).
For the elementary abelian group \(C_p^r\), where \(p\) is prime and \(r\geq 1\), the density is \(p^{1-r}\).
\end{abstract}


\section{Introduction and Preliminaries}\label{sec:framework}
The Cayley graph and the power graph have two different relations on a finite group: multiplication by prescribed generators and containment in cyclic subgroups. Their common vertex set permits comparison without choosing a graph isomorphism or relabelling. Background on power graphs is given in \cite{survey,cameron,cameronghosh,recentsurvey}. Component and quotient structure for power graphs of finite groups, including alternating and symmetric groups, is treated in \cite{bubboloni2017quotient,bubboloni2017alternating}. The undirected power graph of a finite group is studied in \cite{pourgholi2015}, and finite nilpotent groups determined by the enhanced power graph in \cite{mirzargar2026}. Chattopadhyay and Panigrahi \cite{chattopadhyay} study relations between power graphs and Cayley graphs of finite cyclic groups. Explicit enumerative work on cyclic, dihedral, and quaternion-type power graphs includes \cite{hosoya}; other recent directions include power-graph isomorphism and directed reconstruction \cite{das2024,bubboloni2025}.

We ask the following quantitative question:
\begin{quote}
How large a proportion of the edges of a \emph{connected} Cayley graph on a finite group can be retained in its power graph?
\end{quote}
The connectedness requirement is essential. Without it, a single central compatible involution can give proportion one, although the corresponding Cayley graph is merely a perfect matching. Generalized quaternion groups illustrate this degeneracy. Requiring the connection set to generate the group forces a trade-off between elements with favourable local counts and elements needed to connect different parts of the group.

The term ``approximation'' is used in a specified, one-sided sense: the numerator counts retained Cayley edges, not all power-graph edges. We do not identify this objective with minimizing symmetric difference. The distinction, and the fixed-valency relationship between these objectives, is recorded in Remark~\ref{rem:objectives}.

\subsection*{Notation and preliminary results}
All groups and graphs in this paper are finite. We denote $e$ for the identity, $Z(G)$ for the center of $G$, and $\gen{x}$ for the cyclic subgroup generated by $x$. A \emph{connection set} is a subset $S\subseteq G\setminus\{e\}$ satisfying $S=S^{-1}$. Note that possibly $S=\varnothing$.

The simple undirected Cayley graph $\Cay(G,S)$ has vertex set $G$, with
\[
 g\sim h \quad\Longleftrightarrow\quad hg^{-1}\in S.
\]
Thus its edges are $\{g,sg\}$, where $g\in G$ and $s\in S$. This left-multiplication convention is used throughout. The Cayley graph $\Cay(G, S)$ is connected if and only if $\gen{S}=G$.

The undirected power graph $\Pow(G)$ has vertex set $G$. Two distinct vertices $x,y$ are adjacent if $x\in\gen{y}$ or $y\in\gen{x}$. Power adjacency implies commutativity.
The weaker condition that $\gen{x,y}$ is cyclic is not sufficient: elements of orders $2$ and $3$ in $C_6$ generate a cyclic subgroup but are not power-adjacent. That weaker relation is the edge rule of the enhanced power graph; its structure is described in \cite{aalipour2017}, and later surveys and neighbourhood results appear in \cite{ma2022,lewis2024}. Exact connectivity formulas for enhanced power graphs of finite nilpotent groups are given in \cite{bera2025}.
Note that the local count $c_G(s)$ below is a degree in $\Pow(G)$, not in the enhanced power graph.

\begin{definition}
A connection set $S$ is \emph{compatible} if $\Cay(G,S)$ is a spanning subgraph of $\Pow(G)$. The group $G$ is \emph{strongly compatible} if some compatible connection set generates $G$. Define
\[
 \Comp(G)=\{s\in G\setminus\{e\}:g\text{ and }sg\text{ are power-adjacent for every }g\in G\}.
\]
The trivial group is strongly compatible via the empty connection set.
\end{definition}

For $s\ne e$, define
\[
 \Omega_s=\{g\in G:g\in\gen{sg}\text{ or }sg\in\gen{g}\},
 \qquad c_G(s)=|\Omega_s|,
\]
and write
\[
 I(G,S)=|E(\Cay(G,S))\cap E(\Pow(G))|.
\]
\begin{proposition}\label{prop:count}
For every connection set $S$,
\begin{equation}\label{eq:count}
 I(G,S)=\frac12\sum_{s\in S}c_G(s).
\end{equation}
Moreover, $c_G(s)=c_G(s^{-1})$ and $c_G(s)\geq2$ for $s\ne e$.
\end{proposition}
\begin{proof}
The pairs $(g,s)$ with $s\in S$ and $g\in\Omega_s$ represent exactly
the edges in $E(\operatorname{Cay}(G,S))\cap E(P(G))$.
Each edge $\{g,sg\}$ is represented twice, by $(g,s)$ and
$(sg,s^{-1})$, which proves \eqref{eq:count}.

For any $s\ne e$, the map $g\mapsto sg$ is a bijection from
$\Omega_s$ to $\Omega_{s^{-1}}$, so $c_G(s)=c_G(s^{-1})$.
Finally, $e$ and $s^{-1}$ are distinct elements of $\Omega_s$,
since the corresponding edges each join the identity to a
nonidentity element. Thus $c_G(s)\geq 2$.
\end{proof}

For a nontrivial finite group define
\[
 \rho(G)=\max_{\substack{\varnothing\ne S\subseteq G\setminus\{e\}\\S=S^{-1}}}
             \frac{I(G,S)}{|G||S|/2}.
\]
The empty set is excluded to avoid division by zero.

\begin{theorem}\label{thm:density}
For every nontrivial finite group,
\[
 \rho(G)=\max_{s\ne e}\frac{c_G(s)}{|G|}.
\]
Hence $2/|G|\leq\rho(G)\leq1$, and $\rho(G)=1$ if and only if $\Comp(G)\ne\varnothing$.
\end{theorem}
\begin{proof}
By \eqref{eq:count}, the density associated with $S$ is the average of $c_G(s)/|G|$ over $s\in S$,
and cannot exceed the largest term. Conversely, choose a maximizing $s$ and use $S=\{s,s^{-1}\}$, whose local counts agree. The other assertions are immediate.
\end{proof}

\begin{definition}
For a nontrivial finite group, the \emph{connected intersection density} is
\begin{equation}\label{eq:rhogen}
 \rho_{\mathrm{gen}}(G)=
 \max_{\substack{S=S^{-1}\subseteq G\setminus\{e\}\\\gen{S}=G}}
 \frac{I(G,S)}{|G||S|/2}
 =\frac1{|G|}\max_{\substack{S=S^{-1}\subseteq G\setminus\{e\}\\\gen{S}=G}}
 \frac{\sum_{s\in S}c_G(s)}{|S|}.
\end{equation}
The maximum exists because $G\setminus\{e\}$ is an admissible set and there are finitely many connection
sets. The trivial group is excluded from density definitions because its Cayley graph has no edges.
\end{definition}

\begin{proposition}\label{prop:basicconnected}
For every nontrivial finite group,
\[
 \frac2{|G|}\leq\rho_{\mathrm{gen}}(G)\leq\rho(G)\leq1.
\]
If $G=C_p^r$, where $p$ is prime and $r\geq1$, then
\[
 \rho_{\mathrm{gen}}(C_p^r)=\rho(C_p^r)=p^{1-r}.
\]
\end{proposition}

\begin{proof}
By \eqref{eq:count}, the density associated with a nonempty connection
set $S$ is the average of $c_G(s)/|G|$ over all $s\in S$. Since
$2\leq c_G(s)\leq |G|$ for every $s\ne e$, each such density lies
between $2/|G|$ and $1$. Moreover, every generating connection set
is a nonempty connection set, so
$\rho_{\mathrm{gen}}(G)\leq\rho(G)$.

Now let $G=C_p^r$ and fix $s\ne e$. We claim that
$\Omega_s=\gen{s}$. Suppose first that $g\in\Omega_s$. If $g=e$ or
$sg=e$, then $g=e$ or $g=s^{-1}$, respectively, so $g\in\gen{s}$.
Otherwise, $g$ and $sg$ are nonidentity elements. Since they are
power-adjacent, one of the cyclic subgroups $\gen{g}$ and $\gen{sg}$
is contained in the other. Both subgroups have order $p$, so they
are equal. Their common subgroup contains $s=(sg)g^{-1}$ and
therefore equals $\gen{s}$. Hence $g\in\gen{s}$.

Conversely, if $g\in\gen{s}$, then $sg\in\gen{s}$. The two elements
are distinct because $s\ne e$. If one is the identity, it is a
power of the other. If neither is the identity, both generate
$\gen{s}$, so each is a power of the other. Thus $g\in\Omega_s$.

It follows that $c_G(s)=|\gen{s}|=p$ for every $s\ne e$. Hence every
nonempty connection set has density $p/|G|=p^{1-r}$. In particular,
both $\rho_{\mathrm{gen}}(G)$ and $\rho(G)$ equal $p^{1-r}$.
\end{proof}

\begin{proposition}\label{prop:equal}
Let $G$ be a nontrivial finite group and let
\[
 M=\{s\in G\setminus\{e\}:c_G(s)=\max_{x\ne e}c_G(x)\}.
\]
Then $\rho_{\mathrm{gen}}(G)=\rho(G)$ if and only if $\gen{M}=G$.
\end{proposition}
\begin{proof}
By \eqref{eq:count}, the density associated with a nonempty
connection set $S$ is the average of $c_G(s)/|G|$ over $s\in S$.
By Theorem~\ref{thm:density}, the largest possible summand is
$\rho(G)$, and it is attained precisely when $s\in M$. Thus the
average equals $\rho(G)$ if and only if $S\subseteq M$.

By Proposition~\ref{prop:count}, $c_G(s)=c_G(s^{-1})$, so $M$ is
inverse-closed. If $\gen{M}=G$, then $M$ is a generating connection
set whose density is $\rho(G)$. Hence
$\rho_{\mathrm{gen}}(G)=\rho(G)$.

Conversely, suppose that $\rho_{\mathrm{gen}}(G)=\rho(G)$.
Choose a generating connection set $S$ attaining
$\rho_{\mathrm{gen}}(G)$. Since its density is $\rho(G)$, we have
$S\subseteq M$. Therefore $G=\gen{S}\subseteq\gen{M}\subseteq G$,
so $\gen{M}=G$.
\end{proof}


\begin{remark}\label{rem:objectives}
The connected intersection density measures the proportion of Cayley
graph edges that also belong to the power graph. An edge of the power
graph that is absent from the Cayley graph is counted in neither the
numerator nor the denominator. In contrast, the symmetric difference
counts every edge that belongs to exactly one of the two graphs.

Put $n=|G|$ and $t=|E(\Pow(G))|$. For each feasible valency $k$,
that is, each $k$ for which there exists a generating connection set
of size $k$, define
\[
 \lambda_G(k)=
 \max_{\substack{S=S^{-1},\ e\notin S\\ \gen{S}=G,\ |S|=k}}
 \sum_{s\in S}c_G(s).
\]
At fixed valency $k$, the Cayley graph has $nk/2$ edges, and
\eqref{eq:count} shows that the largest possible number of common
edges is $\lambda_G(k)/2$. Therefore
\[
 \rho_{\mathrm{gen}}(G)=\max_k\frac{\lambda_G(k)}{nk},
\]
where the maximum is taken over the feasible valencies.

If the two graphs have $I=I(G,S)$ common edges, then their edge sets
have symmetric difference of size $t+nk/2-2I$. Thus at fixed $k$,
minimizing the size of the symmetric difference is equivalent to
maximizing $I$, and its minimum is
\[
 t+\frac{nk}{2}-\lambda_G(k).
\]

By inclusion--exclusion, the union of the two edge sets has size
$t+nk/2-I$. Hence the proportion of edges in the union that belong
to both graphs is $I/(t+nk/2-I)$. At fixed $k$, this proportion
increases with $I$, since its numerator increases and its positive
denominator decreases. Substituting $I=\lambda_G(k)/2$ gives the
maximum
\[
 \frac{\lambda_G(k)}{nk+2t-\lambda_G(k)}.
\]

When $k$ is allowed to vary, the three objectives need not have the
same optimizing connection sets. The results in this paper concern
the proportion of Cayley graph edges that also belong to the power
graph.
\end{remark}

\section{The density-one boundary: compatible elements}\label{sec:compatible}
In this section we determine which connection elements lie in every power-graph translate, and then characterizes the groups for which a generating connection set retains every Cayley edge. Compatibility forces centrality and prime-power order, together with uniqueness of a subgroup of that prime order. The resulting description of $\Comp(G)$ is the input to the density-one theorem.

\begin{proposition}\label{prop:compatibility-basics}
\begin{enumerate}[label=\textup{(\roman*)},ref=\textup{(\roman*)}]
\item\label{prop:largest} The set $\Comp(G)$ is inverse-closed. A connection set $S$ is compatible if and only if $S\subseteq\Comp(G)$. Thus $\Comp(G)$ is the unique largest compatible connection set.
\item\label{prop:center} For every finite group $G$, one has $\Comp(G)\subseteq Z(G)$.
\item\label{lem:primepower} Every element of $\Comp(G)$ has prime-power order.
\item\label{lem:unique} If $s\in\Comp(G)$ has order $p^k$, then $G$ has a unique subgroup of order $p$.
\item\label{lem:abelianp} A nontrivial finite abelian $p$-group with a unique subgroup of order $p$ is cyclic.
\end{enumerate}
\end{proposition}
\begin{proof}
\begin{enumerate}[label=\textup{(\roman*)},ref=\textup{(\roman*)}]
\item If $s\in\Comp(G)$, apply the defining condition to $s^{-1}g$. The resulting adjacent pair is $\{s^{-1}g,g\}$, proving $s^{-1}\in\Comp(G)$. The other assertions follow by considering the edges associated with each individual connection element.
\item For $s\in\Comp(G)$ and $g\in G$, the vertices $g$ and $sg$ commute. Hence $gsg=sg^2$, and cancellation gives $gs=sg$.
\item Let $s\in\Comp(G)$ and restrict the compatibility condition to $H=\gen{s}$. Suppose that distinct primes $p,q$ divide $|s|$. Decompose the cyclic group $H$ into its Sylow factors and write $s=(s_r)_r$. Choose $g\in H$ with $p$-component $e$, $q$-component $s_q^{-1}$, and all other components equal to $e$. The $p$-component of $sg$ is nontrivial, whereas that of $g$ is trivial, so $sg\notin\gen{g}$. The $q$-component of $g$ is nontrivial, whereas that of $sg$ is trivial, so $g\notin\gen{sg}$. This contradicts compatibility.
\item By part~\ref{prop:center}, $s$ is central. Let $t$ have order $p$. Suppose $t\notin\gen{s}$. Then $\gen{t}\cap\gen{s}=\{e\}$. The inclusion $st\in\gen{t}$ would imply $s\in\gen{t}$ and hence $t\in\gen{s}$, a contradiction. On the other hand, if $t=(st)^j$, centrality gives $s^j=t^{1-j}$. Trivial intersection forces $p^k\mid j$ and $p\mid(1-j)$, which is impossible. Thus $t$ and $st$ are not power-adjacent, contrary to compatibility. Every element of order $p$ therefore belongs to $\gen{s}$, which has a unique subgroup of order $p$.
\item By the decomposition theorem for finite abelian groups, a noncyclic abelian $p$-group has at least two nontrivial cyclic direct factors. Their subgroups of order $p$ are distinct.
\end{enumerate}
\end{proof}

For a prime $p\mid |G|$, write $O_p(Z(G))$ for the subgroup of all $p$-elements of $Z(G)$, including $e$. Since the center is abelian, this is its Sylow $p$-subgroup, and also its largest normal $p$-subgroup. Let $\mathcal{U}(G)$ be the set of primes $p$ for which $G$ has a unique subgroup of order $p$.

\begin{theorem}\label{thm:main}
For $s\ne e$, the following are equivalent:
\begin{enumerate}[label=\textup{(\roman*)}]
\item $s\in\Comp(G)$;
\item $s\in Z(G)$, $|s|=p^k$ for a prime $p$ and $k\geq1$, and $G$ has a unique subgroup of order $p$.
\end{enumerate}
Consequently,
\begin{equation}\label{eq:main}
 \Comp(G)=\bigcup_{p\in \mathcal{U}(G)}\bigl(O_p(Z(G))\setminus\{e\}\bigr),
 \qquad
 |\Comp(G)|=\sum_{p\in \mathcal{U}(G)}\bigl(|O_p(Z(G))|-1\bigr).
\end{equation}
\end{theorem}
\begin{proof}
Necessity follows from Proposition~\ref{prop:compatibility-basics}, parts~\ref{prop:center}--\ref{lem:unique}.

For sufficiency, suppose $s$ satisfies (ii). Fix $g\in G$. Within $\gen{g}$ write $g=g_p g_{p'}$, where $g_p$ has $p$-power order, $g_{p'}$ has order coprime to $p$, and the factors commute. The subgroup $A=\gen{s,g_p}$ is an abelian $p$-group. It has a unique subgroup of order $p$, since $G$ does. By Proposition~\ref{prop:compatibility-basics}, part~\ref{lem:abelianp}, $A$ is cyclic. Its subgroups are linearly ordered, so $\gen{g_p}$ and $\gen{sg_p}$ are comparable.

The element $g_{p'}$ commutes with $A$ and has order coprime to $|A|$. Therefore
\[
 \gen{g}=\gen{g_p}\times\gen{g_{p'}},\qquad
 \gen{sg}=\gen{sg_p}\times\gen{g_{p'}}.
\]
These cyclic subgroups are comparable, so $g$ and $sg$ are power-adjacent. Since $g$ was arbitrary, $s\in\Comp(G)$. Formula~\eqref{eq:main} follows, since nonidentity elements of different prime-power orders lie in disjoint parts of the union.
\end{proof}

\begin{corollary}\label{cor:strong}
A finite group is strongly compatible if and only if it is cyclic, including the trivial group.
\end{corollary}
\begin{proof}
Suppose a compatible set $S$ generates $G$. By Proposition~\ref{prop:compatibility-basics}, part~\ref{prop:center}, $G=\gen{S}\subseteq Z(G)$, so $G$ is abelian. In an abelian group, uniqueness of a subgroup of order $p$ is equivalent to cyclicity of the Sylow $p$-subgroup, by Proposition~\ref{prop:compatibility-basics}, part~\ref{lem:abelianp}. Theorem~\ref{thm:main} says that every element of $S$ is supported only in a cyclic Sylow factor. If a Sylow factor were noncyclic, the projection of $\gen{S}$ to that factor would be trivial. Thus every Sylow factor is cyclic, and $G$ is cyclic.

Conversely, if $G$ is cyclic, Theorem~\ref{thm:main} shows that the union of the nonidentity elements of its Sylow subgroups is compatible. This union generates $G$. For the trivial group use the empty set.
\end{proof}

\begin{remark}
The prime-power restriction applies to each individual connection element, not to the order of a strongly compatible group. For $C_6=\gen{a}$, the set $\{a^2,a^3,a^4\}$ is compatible and generates $C_6$.
\end{remark}

\begin{corollary}\label{cor:abelian}
Let $G$ be abelian and let $S_p$ denote its Sylow $p$-subgroup for each prime $p\mid |G|$. In the decomposition $G=\prod_p S_p$, identify each factor with its natural embedded subgroup.
\[
 \Comp(G)=\bigcup_{\substack{p\mid |G|\\S_p\text{ cyclic}}}(S_p\setminus\{e\}),
 \qquad
 |\Comp(G)|=\sum_{\substack{p\mid |G|\\S_p\text{ cyclic}}}(|S_p|-1).
\]
In particular, if $n=\prod_{p\mid n}p^{a_p}$, then
\[
 |\Comp(C_n)|=\sum_{p\mid n}(p^{a_p}-1).
\]
\end{corollary}
\begin{proof}
Apply Proposition~\ref{prop:compatibility-basics}, part~\ref{lem:abelianp} and Theorem~\ref{thm:main}.
\end{proof}
Thus every nonidentity element of a cyclic $p$-group is compatible, whereas a noncyclic abelian $p$-group has no nonempty compatible connection set.

\begin{corollary}\label{cor:nilpotent}
Let $G$ be nilpotent and write $G=\prod_p S_p$, where $S_p$ is the Sylow $p$-subgroup of $G$. Then
\[
 \Comp(G)=\bigcup_{\substack{p\mid |G|\\S_p\text{ has a unique subgroup of order }p}}
                  (Z(S_p)\setminus\{e\}),
\]
where each center is embedded in its corresponding direct factor.
\end{corollary}
\begin{proof}
Every element of order $p$ lies entirely in $S_p$, and $Z(G)=\prod_p Z(S_p)$. Apply Theorem~\ref{thm:main}.
\end{proof}

\begin{example}
The quaternion group $Q_8$ has center $\{1,-1\}$ and a unique subgroup of order $2$. Thus $\Comp(Q_8)=\{-1\}$, which is compatible but does not generate $Q_8$. More generally, a generalized quaternion $2$-group has compatible set consisting of its unique involution. Here ``generalized quaternion'' uses the convention that the order is a power of $2$; the statement is not made for every dicyclic group under the notation $Q_{4n}$.
\end{example}

\begin{example}
If $Z(G)=\{e\}$, then $\Comp(G)=\varnothing$. Examples include $S_3$, $A_4$, and the dihedral group of order $10$.
\end{example}

\begin{example}
For $G=C_6\times C_2$, the Sylow $2$-subgroup is $C_2\times C_2$, while the Sylow $3$-subgroup is cyclic. Therefore $\Comp(G)$ consists of the two nonidentity elements of the embedded $C_3$ factor. A noncyclic abelian group can consequently have a nonempty compatible connection set.
\end{example}

\begin{theorem}\label{thm:connectedone}
For a nontrivial finite group $G$,
\[
 \rho_{\mathrm{gen}}(G)=1\quad\Longleftrightarrow\quad G\text{ is cyclic}.
\]
\end{theorem}
\begin{proof}
An admissible set has density one exactly when all its Cayley edges belong to the power graph. Thus density one is equivalent to strong compatibility, and Corollary~\ref{cor:strong} applies.
\end{proof}

The distinction from unrestricted density is already visible in $Q_8$: the set consisting of its unique involution has density one but does not generate the group. Theorem~\ref{thm:quaternionconnected} gives the connected value $3/5$.

\begin{corollary}\label{thm:amplification}
Let $T$ be a generating connection set of a noncyclic finite group $H$.
If $m=q^a$, where $q$ is prime, $q\nmid |H|$, and $a\geq1$, then
\[
 \rho_{\mathrm{gen}}(H\times C_m)\geq \frac{m-1}{m-1+|T|}.
\]
Thus noncyclic groups can have connected intersection density arbitrarily close to $1$.
\end{corollary}
\begin{proof}
By Theorem~\ref{thm:main}, every $(e,z)$ with $z\ne e$ is compatible:
it is central, has prime-power order, and $H\times C_m$ has a unique
subgroup of order $q$. These $m-1$ elements, together with $T\times\{e\}$,
form a generating connection set and give the bound. The product is
noncyclic because it projects onto $H$, and the bound tends to $1$ as
$a\to\infty$.
\end{proof}

\section{Generating constraints and level-set optimization}\label{sec:optimization}
Dihedral and dicyclic calculations maximize an average of local counts over inverse-closed generating sets. This section separates that optimization from the group structure. The first subsection treats weights that attain their maximum on a subgroup and are constant outside it. The second reformulates the problem on inversion orbits and gives a parametric certificate for the optimal value.

\subsection{A general level-set principle}
The statement concerns arbitrary inverse-symmetric weights, so its hypotheses and its dependence on generation are explicit.

\begin{theorem}\label{thm:threshold}
Let $G$ be nontrivial and let $w:G\setminus\{e\}\to\mathbb R$ satisfy $w(s)=w(s^{-1})$. Let $H<G$ be a nontrivial subgroup. Suppose that
\begin{enumerate}[label=\textup{(\roman*)}]
\item $M=\max_{s\ne e}w(s)$ and the set
\[
 B=\{h\in H\setminus\{e\}:w(h)=M\}
\]
generates $H$;
\item $w(s)=d<M$ for every $s\in G\setminus H$.
\end{enumerate}
Define
\[
 q=\min\{|U|:U=U^{-1}\subseteq G\setminus H,\ \gen{H,U}=G\}.
\]
For each distinct weight $t$ on $H\setminus\{e\}$, put
\[
 T_t=\{h\in H\setminus\{e\}:w(h)\geq t\}.
\]
Then
\begin{equation}\label{eq:threshold}
 \max_{\substack{S=S^{-1},\ e\notin S\\\gen{S}=G}}
 \frac{\sum_{s\in S}w(s)}{|S|}
 =\max_t\frac{\sum_{h\in T_t}w(h)+qd}{|T_t|+q}.
\end{equation}
For a maximizing level $t$, adjoining any $q$-element set $U$ from the definition of $q$ gives an optimal connection set $T_t\cup U$.
\end{theorem}
\begin{proof}
The set $B$ is inverse-closed. Choosing $B$ together with a minimum set $U$ gives an admissible average strictly greater than $d$. Every generating set meets $G\setminus H$, so every admissible average is strictly less than $M$.

Let $S$ be optimal and let its average be $a$. Thus $d<a<M$. Every missing element of $B$ can be added together with its inverse; this would strictly increase the average, a contradiction. Hence $B\subseteq S$, so $S\cap H$ generates $H$.

Write $R=S\cap H$ and $V=S\setminus H$. Then $\gen{H,V}=G$, so $|V|\geq q$. Replacing $V$ by a minimum set $U$ preserves generation. Since all outside weights equal $d<a$, reducing their number strictly increases the average unless $|V|=q$. An optimal set therefore has exactly $q$ outside elements.

With these $q$ outside elements fixed, every omitted inside element of weight greater than $a$ could be added, and every included inside element of weight less than $a$ could be removed, strictly improving the average. Removal preserves generation because all elements of $B$ remain. Inverse symmetry permits these changes in whole inverse orbits. Elements of weight $a$ may all be included without changing the average. Consequently an optimal inside set can be chosen as $T_t$ for a weight level $t$. Conversely, every $T_t$ contains $B$ and generates $H$, so $T_t\cup U$ is admissible. This proves the formula and the attainment assertion.
\end{proof}

\begin{corollary}\label{cor:subgroupbound}
Under the hypotheses of Theorem~\ref{thm:threshold}, writing $h=|H|$, the optimal average is at most
\[
 \frac{M(h-1)+dq}{h-1+q}.
\]
If every nonidentity element of $H$ has weight $M$, equality holds.
\end{corollary}
\begin{proof}
Each superlevel set has at most $h-1$ elements and weight sum at most $M|T_t|$. Since $M>d$, the expression $(Mr+dq)/(r+q)$ increases with $r$. If all inside weights equal $M$, take $T=H\setminus\{e\}$.
\end{proof}

For $w=c_G$, the optimal average in Theorem~\ref{thm:threshold} is $|G|\rho_{\mathrm{gen}}(G)$. The theorem is not asserted for an arbitrary subgroup. Its hypotheses, generation by the maximum-weight elements and constant weight outside the subgroup, hold for the rotation subgroup of every dihedral group $D_{2n}$, $n\ge 3$; this is why \eqref{eq:dihedralgeneral} is not restricted to prime powers. They fail for the cyclic index-two subgroup of a generalized quaternion group, whose unique maximum-weight element generates only the order-two subgroup; Section~\ref{sec:dic} supplies the argument used there. They also fail for groups such as $C_4\times C_2$. The values for those groups in Section~\ref{sec:computation} come from exact dynamic programming, not from Theorem~\ref{thm:threshold}.

\subsection{An exact inverse-orbit formulation}
Let $O_1,\ldots,O_m$ be the orbits of inversion on $G\setminus\{e\}$. Put
\[
 a_i=|O_i|,\qquad b_i=\sum_{s\in O_i}c_G(s).
\]
Every connection set is a union of these orbits. Thus
\begin{equation}\label{eq:orbitopt}
 |G|\rho_{\mathrm{gen}}(G)=
 \max_{\substack{J\subseteq\{1,\ldots,m\}\\\gen{\bigcup_{j\in J}O_j}=G}}
 \frac{\sum_{j\in J}b_j}{\sum_{j\in J}a_j}.
\end{equation}
The generation condition is not discarded in this formula. Unlike unrestricted density, the connected invariant need not be attained by one inversion orbit.

\begin{proposition}\label{prop:parametric}
For $t\in\mathbb R$, let $K_t$ be the subgroup generated by all orbits with $b_i-ta_i\geq0$. Define
\[
 P(t)=\sum_i\max\{b_i-ta_i,0\},
\]
and let $C(t)$ be the minimum of
\[
 \sum_{i\in J}(ta_i-b_i)
\]
over sets $J$ of negative-weight orbits ($b_i-ta_i<0$) such that $K_t$ together with these orbits generates $G$. The empty set is allowed when $K_t=G$. Then
\[
 P(t)-C(t)=
 \max_{\substack{S=S^{-1},\ e\notin S\\\gen{S}=G}}
 \sum_{s\in S}(c_G(s)-t).
\]
This expression is positive, zero, or negative according as $t$ is less than, equal to, or greater than $|G|\rho_{\mathrm{gen}}(G)$.
\end{proposition}
\begin{proof}
Adding every positive-weight orbit increases the objective and cannot destroy generation; adding all zero-weight orbits has no effect on the objective. After doing so, the only remaining choice is a generating extension of minimum negative penalty, which gives $P(t)-C(t)$. For a nonempty admissible $S$, the objective is $|S|$ times its average local count minus $t|S|$. Comparing $t$ with the maximum admissible average proves the sign statement.
\end{proof}

This certificate separates favourable orbits from the cost of completing generation. Computing the minimum generating extension can still require exponential search; no polynomial-time algorithm for arbitrary groups is claimed.

\section{Connected approximation in dihedral groups}\label{sec:dihedral}
In a dihedral group the rotation subgroup carries the cyclic local counts computed in the first subsection, while every reflection has local count $2$. Connectedness forces at least one reflection into the connection set. This section records those local counts and then evaluates the connected density by the level-set principle. When the rotation order is a prime power, the maximum uses every nonidentity rotation together with one reflection.

\subsection{Local weights in cyclic groups}\label{sec:cyclic}
Local counts on cyclic groups are the numerical input to the dihedral and dicyclic formulas. This subsection computes them from the Sylow decomposition of the element. An element and its translate are power-adjacent precisely when no primary component strictly decreases in order, and the resulting product formula also recovers which cyclic elements are compatible.

Write $C_n$ additively as $\mathbb Z/n\mathbb Z$, where $n=\prod_{i=1}^r p_i^{a_i}$. For $s\in C_n$, suppose its component in $\mathbb Z/p_i^{a_i}\mathbb Z$ has order $p_i^{b_i}$, where $0\leq b_i\leq a_i$. Put
\[
 N_i=p_i^{a_i},\qquad
 u_i=\begin{cases}0,&b_i=0,\\p_i^{b_i-1},&b_i>0.\end{cases}
\]
\begin{theorem}\label{thm:cycliccount}
For $s\ne0$,
\begin{equation}\label{eq:cycliccount}
 c_{C_n}(s)=2\prod_{i=1}^r(N_i-u_i)-\prod_{i=1}^r(N_i-2u_i).
\end{equation}
Consequently,
\[
 I(C_n,S)=\frac12\sum_{s\in S}
 \left[2\prod_{i=1}^r(N_i-u_i(s))-\prod_{i=1}^r(N_i-2u_i(s))\right].
\]
\end{theorem}
\begin{proof}
In a cyclic $p$-group, subgroups are linearly ordered, so containment is equivalent to comparison of element orders. If $b_i=0$, adding $s_i$ changes no orders. Otherwise let $t_i=a_i-b_i$, the $p_i$-adic valuation of $s_i$. The order of $x_i+s_i$ is strictly larger than that of $x_i$ exactly when $x_i$ is divisible by $p_i^{t_i+1}$. There are $p_i^{b_i-1}=u_i$ such choices. The involution $x_i\mapsto-s_i-x_i$ exchanges strict increases and strict decreases. There are therefore also $u_i$ strict decreases, with the remaining $N_i-2u_i$ choices preserving order.

Globally, $\gen{x}\subseteq\gen{x+s}$ exactly when no component exhibits a strict decrease. This gives $\prod_i(N_i-u_i)$ choices. Reverse containment has the same count. Both hold exactly when every component preserves order, giving $\prod_i(N_i-2u_i)$ choices. Inclusion--exclusion proves \eqref{eq:cycliccount}; Proposition~\ref{prop:count} gives the intersection formula.
\end{proof}

If $s$ is supported on exactly one Sylow component, \eqref{eq:cycliccount} gives $c(s)=n$. If at least two components are nonzero, there are choices with a strict increase at one prime and a strict decrease at another; hence $c(s)<n$. This independently confirms the cyclic case of Theorem~\ref{thm:main}.

We use $D_{2n}$ for the dihedral group of order $2n$, with $n\geq3$:
\[
 D_{2n}=\gen{a,b\mid a^n=b^2=e,\ bab^{-1}=a^{-1}}.
\]
For cyclic local counts below, $c_{C_m}(a^k)$ refers to the corresponding element of the rotation subgroup, and is explicitly given by Theorem~\ref{thm:cycliccount}. The power graphs of these families have also been studied from other enumerative viewpoints; see \cite{hosoya}.

\begin{proposition}\label{prop:dihedral}
For $D_{2n}$ with $n\geq3$, one has
\[
 \Comp(D_{2n})=\varnothing,\qquad \rho(D_{2n})=\frac12.
\]
More precisely, for $1\leq k<n$ and $0\leq j<n$,
\[
 c_{D_{2n}}(a^k)=c_{C_n}(a^k),\qquad
 c_{D_{2n}}(a^jb)=2.
\]
Thus every inverse-closed connection set satisfies
\[
 I(D_{2n},S)=\frac12\left(
 \sum_{a^k\in S\cap\gen{a}}c_{C_n}(a^k)
 +2|S\setminus\gen{a}|\right).
\]
\end{proposition}
\begin{proof}
For a nonidentity rotation $s$, rotation vertices contribute exactly the cyclic local count. For a reflection $g$, the vertices $g$ and $sg$ are distinct involutions, hence are not power-adjacent. For a reflection $s$, multiplication interchanges the two cosets of $\gen{a}$. A rotation and a reflection are power-adjacent only if the rotation is $e$. The only successful choices are consequently $g=e,s$.

The cyclic local count is at most $n$, and equals $n$ for any nonidentity rotation of prime-power order by Corollary~\ref{cor:abelian}. Such a rotation exists. Therefore the maximum local count is $n$, giving density $1/2$ and an empty compatible set. The intersection identity follows from Proposition~\ref{prop:count}.
\end{proof}

\begin{theorem}\label{thm:dihedralconnected}
For $n\geq3$, put $v_k=c_{C_n}(a^k)$, $1\leq k<n$, as given by Theorem~\ref{thm:cycliccount}. For each distinct value $t$ among the $v_k$, let
\[
 R_t=\{a^k:1\leq k<n,\ v_k\geq t\}.
\]
Then
\begin{equation}\label{eq:dihedralgeneral}
 \rho_{\mathrm{gen}}(D_{2n})=
 \max_t\frac{2+\sum_{a^k\in R_t}v_k}{2n(1+|R_t|)}.
\end{equation}
Each $R_t$ generates $\gen{a}$, and an optimal connected Cayley graph is obtained from $S=R_t\cup\{b\}$ at a maximizing level.

If $n$ is a prime power, this simplifies to
\begin{equation}\label{eq:dihedralprimepower}
 \rho_{\mathrm{gen}}(D_{2n})=\frac{n(n-1)+2}{2n^2}.
\end{equation}
Equation \eqref{eq:dihedralprimepower} holds if and only if $n$ is a prime power. In all cases $\rho_{\mathrm{gen}}(D_{2n})<1/2=\rho(D_{2n})$.
\end{theorem}
\begin{proof}
Take $H=\gen{a}$ and $w=c_{D_{2n}}$. Proposition~\ref{prop:dihedral} gives outside weight $d=2$ and inside maximum $M=n$. Every nonidentity element of each Sylow subgroup of $H$ has cyclic local count $n$. Their union generates the cyclic group $H$. A single reflection is inverse-closed and, together with $H$, generates $D_{2n}$, so $q=1$. Theorem~\ref{thm:threshold} proves \eqref{eq:dihedralgeneral} and the attainment statement.

For prime-power $n$, every nonidentity rotation has local count $n$, so Corollary~\ref{cor:subgroupbound} gives \eqref{eq:dihedralprimepower}. Conversely, the sentence after Theorem~\ref{thm:cycliccount} says that $v_k=n$ for every $1\le k<n$ precisely when $n$ has one prime factor. For $n>2$ the quotient $(nr+2)/(r+1)$ is strictly increasing in $r$. A level containing $r$ rotations therefore has average at most $(n(n-1)+2)/n$, with equality only for $r=n-1$ and every weight equal to $n$. Division by $2n$ gives the criterion. Finally, every generating connection set contains a reflection of local count $2<n$, while no local count exceeds $n$, so the average is strictly less than $n$.
\end{proof}

\begin{example}
For $D_{12}$, the nonidentity rotations in the Sylow subgroups of $C_6$ have local count $6$ (three elements), while its two order-six generators have count $4$. The optimal set consists of those three high-count rotations and one reflection. Thus
\[
 \rho_{\mathrm{gen}}(D_{12})=\frac{3\cdot6+2}{12\cdot4}=\frac5{12}.
\]
Using every nonidentity rotation and one reflection instead gives $7/18$, which is smaller. Therefore the prime-power formula cannot be extended to arbitrary $n$ merely by including all rotations.
\end{example}

\section{Dicyclic and generalized quaternion groups}\label{sec:dic}
Dicyclic groups have a unique involution, and that involution is compatible, so the unrestricted density equals $1$. A generating connection set must nevertheless leave the cyclic subgroup of index two. This section computes the local counts and then optimizes the connected density. In the generalized quaternion subfamily every nonidentity rotation has the same cyclic count, and the level-set maximum collapses to a closed formula.

For $n\geq2$, define the dicyclic group of order $4n$ by
\[
 \operatorname{Dic}_n=\gen{a,b\mid a^{2n}=e,\ b^2=a^n,\ bab^{-1}=a^{-1}},
 \qquad z=a^n.
\]
Its unique involution is $z$, and its center is $\{e,z\}$. The notation $Q_{4n}$ is sometimes used for this entire family. Here we reserve ``generalized quaternion'' for the subfamily of order a power of $2$, namely $n=2^{m-2}$ and $Q_{2^m}=\operatorname{Dic}_{2^{m-2}}$, where $m\geq3$.

\begin{proposition}\label{prop:dicyclic}
For every $n\geq2$,
\[
 \Comp(\operatorname{Dic}_n)=\{z\},\qquad
 \rho(\operatorname{Dic}_n)=1.
\]
For nonidentity rotations and all elements outside the rotation subgroup,
\[
 c_{\operatorname{Dic}_n}(a^k)=
 c_{C_{2n}}(a^k)+2n\,\mathbf{1}_{\{k=n\}},
 \qquad
 c_{\operatorname{Dic}_n}(a^jb)=4.
\]
In particular,
\[
 I(\operatorname{Dic}_n,S)=\frac12\left(
 \sum_{a^k\in S\cap\gen{a}}c_{C_{2n}}(a^k)
 +2n\,\mathbf{1}_{\{z\in S\}}+4|S\setminus\gen{a}|\right).
\]
\end{proposition}
\begin{proof}
Every element outside $\gen{a}$ has order $4$, square $z$, and cyclic subgroup $\{e,z,g,g^{-1}\}$. Hence $z$ is the unique involution; the presentation also gives $Z(\operatorname{Dic}_n)=\{e,z\}$. Theorem~\ref{thm:main} gives the compatible set and Theorem~\ref{thm:density} gives density one.

For a rotation $s\ne e$, vertices inside $\gen{a}$ contribute the cyclic count. Outside it, $g$ and $sg$ both have order $4$, and are adjacent precisely when $sg=g^{-1}=zg$, or equivalently $s=z$. This adds $2n$ exactly in that case.

For $s=a^jb$, if $g$ is a rotation, $sg$ lies outside $\gen{a}$, and its cyclic subgroup has rotation part $\{e,z\}$. Thus precisely $g=e,z$ contribute. If $g$ lies outside $\gen{a}$, the product $sg$ is a rotation. Power adjacency requires $sg=e$ or $z$, giving precisely $g=s^{-1},s$. These four choices are distinct. The intersection formula follows by summation.
\end{proof}

\begin{corollary}\label{cor:quaternioncounts}
For the generalized quaternion group $Q_{2^m}$, $m\geq3$, let $A=\gen{a}$ have order $2^{m-1}$ and let $z$ be its unique involution. Then
\[
 c(z)=2^m,\qquad
 c(a^k)=2^{m-1}\quad(a^k\notin\{e,z\}),\qquad
 c(a^jb)=4.
\]
The unique largest compatible set is $\{z\}$, and its Cayley graph is a perfect matching contained in the power graph, with $2^{m-1}$ edges. These groups have density one but are not strongly compatible.
\end{corollary}
\begin{proof}
Every nonidentity element of the cyclic $2$-group $A$ has cyclic local count $|A|$. Apply Proposition~\ref{prop:dicyclic}. Multiplication by the involution $z$ partitions the vertices into disjoint pairs. Noncyclicity and Corollary~\ref{cor:strong} rule out strong compatibility.
\end{proof}

\begin{theorem}\label{thm:dicconnected}
Let $G=\operatorname{Dic}_n$, $A=\gen{a}$, and $N=2n$. For $n\geq4$, put
\[
 v_k=c_{C_N}(a^k),\qquad
 R_t=\{a^k:1\leq k<N,\ v_k\geq t\},
\]
where $t$ ranges over the distinct values of $v_k$. Then
\begin{equation}\label{eq:dicgeneral}
 \rho_{\mathrm{gen}}(\operatorname{Dic}_n)=
 \max_t\frac{N+8+\sum_{a^k\in R_t}v_k}{2N(|R_t|+2)}.
\end{equation}
A maximizing level gives an optimal set
\[
 S=R_t\cup\{b,b^{-1}\}.
\]
The two smaller cases are
\[
 \rho_{\mathrm{gen}}(\operatorname{Dic}_2)=\frac35,
 \qquad
 \rho_{\mathrm{gen}}(\operatorname{Dic}_3)=\frac8{15}.
\]
\end{theorem}
\begin{proof}
First suppose $N\geq8$. Since $z$ has local count $2N$ and $G$ is noncyclic, every optimal set contains $z$: otherwise adjoining this involution would strictly increase the average local count.

Write $r=|S\cap A|$ and $u=|S\setminus A|$ for an optimal set. Inverse closure forces $u\geq2$, since every outside element has order four and its inverse is distinct and outside $A$. Proposition~\ref{prop:dicyclic} gives
\begin{equation}\label{eq:dicupper}
 \sum_{s\in S}c_G(s)\leq N(r+1)+4u.
\end{equation}
Consequently the average local count is at most $N$, because
\[
 N(r+1)+4u-N(r+u)=N-u(N-4)\leq8-N\leq0.
\]
All nonidentity elements in the Sylow subgroups of the cyclic group $A$ have cyclic count $N$. Adding all such elements to an optimal set therefore cannot reduce its average, and their union generates $A$. We may hence choose an optimal set whose rotation part generates $A$ and contains all cyclic-count-$N$ elements, including $z$.

The optimal average is strictly greater than $4$: take these rotations together with one outside inverse pair, noting that their actual local counts are at least $N\geq8$. Once the rotations generate $A$, replacing the outside part by $\{b,b^{-1}\}$ preserves generation. Since every outside weight is $4$, this cannot reduce the average and is strictly improving if more than two outside elements were present.

It remains to optimize the rotation part. Keep all cyclic-count-$N$ elements, which already generate $A$, and the fixed outside pair. Including any omitted rotation whose actual local count is above the optimal average, or removing any included lower-count rotation, improves the average. These operations preserve inverse closure and generation. For rotations other than $z$, actual and cyclic local counts agree. The element $z$ is always retained, and its cyclic count $N$ places it in every superlevel set. Ties can all be included without changing the average. Thus an optimal rotation part is some $R_t$. Its actual weight sum is $N+\sum_{R_t}v_k$, and the outside pair contributes $8$, proving \eqref{eq:dicgeneral}. Conversely, each $R_t$ contains all cyclic-count-$N$ elements and generates $A$, so every displayed candidate is admissible.

For $n=2$, every element other than $e,z$ has local count $4$, while $c(z)=8$. An optimal set contains $z$ and must contain at least two distinct inverse pairs of noninvolutions to generate $Q_8$. Thus its size is at least five, and its average local count is at most $4+4/5$, attained by $\{z,a,a^{-1},b,b^{-1}\}$. Dividing by eight gives $3/5$.

For $n=3$, keep $z$, whose count is $12$. The two order-three rotations have count $6$, the two order-six rotations have count $4$, and all six outside elements have count $4$. If the set contains only one outside inverse pair, generation requires at least one nontrivial rotation inverse pair in addition to $z$. Hence its size is at least five and its average is at most $6+(6-4)/5=32/5$: here the bonus of $z$ above baseline $6$ is $6$, and the two outside elements have total deficit $4$. This bound is attained by $z$, the order-three pair, and $\{b,b^{-1}\}$. If there are at least four outside elements, their deficit from baseline $6$ is at least eight, greater than the bonus six, so the average is below six. The optimum is therefore $(32/5)/12=8/15$.
\end{proof}

\begin{theorem}\label{thm:quaternionconnected}
For $m\geq3$, set $N=2^{m-1}$, so $|Q_{2^m}|=2N$. Then
\begin{equation}\label{eq:qconnected}
 \rho_{\mathrm{gen}}(Q_{2^m})=\frac{N^2+8}{2N(N+1)}.
\end{equation}
The set $(A\setminus\{e\})\cup\{b,b^{-1}\}$ is optimal, where $A=\gen{a}$ is the cyclic subgroup of order $N$. Moreover,
\[
 \rho_{\mathrm{gen}}(Q_8)=\frac35,\qquad
 \rho_{\mathrm{gen}}(Q_{16})=\frac12,\qquad
 \rho_{\mathrm{gen}}(Q_{2^m})<\frac12\quad(m\geq5),
\]
and $\rho_{\mathrm{gen}}(Q_{2^m})\to1/2$ as $m\to\infty$, whereas $\rho(Q_{2^m})=1$ for every $m\geq3$.
\end{theorem}
\begin{proof}
For $N\geq8$, every nonidentity element of the cyclic $2$-group $A$ has cyclic local count $N$. Thus the level-set formula in Theorem~\ref{thm:dicconnected} has a single candidate, with $|R|=N-1$ and cyclic weight sum $N(N-1)$. It gives
\[
 \frac{N+8+N(N-1)}{2N(N+1)}=\frac{N^2+8}{2N(N+1)}.
\]
For $N=4$, the small case in that theorem gives the same formula and the stated optimal set. Finally,
\[
 \rho_{\mathrm{gen}}(Q_{2^m})-\frac12
 =\frac{8-N}{2N(N+1)},
\]
which proves all comparisons and the limit. The unrestricted density follows from Proposition~\ref{prop:dicyclic}.
\end{proof}

\begin{remark}
The value at $Q_{16}$ is attained by more than one valency: $\{z,a,a^{-1},b,b^{-1}\}$ has density $1/2$, as does the full nonidentity rotation set together with an outside inverse pair. The theorem asserts an optimal set, not uniqueness of all optimizers.
\end{remark}

\section{The groups $A_4$, $S_4$, and $A_5$}\label{sec:alternating}
The cyclic subgroups of these three groups are cyclic of prime-power order, and each nonidentity element lies in exactly one maximal cyclic subgroup. The local count is therefore the order of that subgroup. The resulting upper bound is attained by an inverse-closed connection set of size $4$.

\begin{lemma}\label{lem:uniquecyclic}
Let $s\in G\setminus\{e\}$, and suppose that $s$ lies in a unique maximal cyclic subgroup $M$. If $M$ is a cyclic $p$-group, then $\Omega_s(G)=M$ and $c_G(s)=|M|$.
\end{lemma}

\begin{proof}
The identity lies in $\Omega_s(G)$. If $g\in\Omega_s(G)\setminus\{e\}$, then $g\in\gen{sg}$ or $sg\in\gen{g}$. In the first case $g^{-1}\in\gen{sg}$, so $s=(sg)g^{-1}\in\gen{sg}$. In the second case $s=(sg)g^{-1}\in\gen{g}$. Thus $s$ lies in a nontrivial cyclic subgroup containing $g$. That subgroup is contained in the unique maximal cyclic subgroup $M$, so $g\in M$.

Conversely, the subgroups of a cyclic $p$-group are totally ordered by inclusion. For any $x,y\in M$, either $x\in\gen{y}$ or $y\in\gen{x}$. Taking $x=s$ and $y=g$ shows that every element of $M$ lies in $\Omega_s(G)$.
\end{proof}

\begin{theorem}\label{thm:alternating}
The connected densities of $A_4$, $S_4$, and $A_5$ are
\[
\rho_{\mathrm{gen}}(A_4)=\frac14,\qquad
\rho_{\mathrm{gen}}(S_4)=\frac16,\qquad
\rho_{\mathrm{gen}}(A_5)=\frac1{12}.
\]
Each value is attained by an inverse-closed generating connection set of size $4$.
\end{theorem}

\begin{proof}
In $A_4$ the nonidentity elements have orders $2$ and $3$. An involution generates a maximal cyclic subgroup of order $2$, while a $3$-cycle generates one of order $3$. Lemma~\ref{lem:uniquecyclic} gives local counts $2$ and $3$, so every connection set has average weight at most $3$. Hence $\rho_{\mathrm{gen}}(A_4)\le 3/12=1/4$, with equality only for an inverse-closed generating set of $3$-cycles. The set
\[
\{(1\,2\,3),\ (1\,3\,2),\ (1\,2\,4),\ (1\,4\,2)\}
\]
is such a set once its generated subgroup is shown to be $A_4$. A subgroup of order $6$ would have index $2$ and would contain every square. Every $3$-cycle is a square, because $(a\,c\,b)^2=(a\,b\,c)$, but $A_4$ has eight $3$-cycles. Thus $A_4$ has no subgroup of order $6$. The two displayed cyclic subgroups intersect trivially and are distinct, so the subgroup they generate has order divisible by $3$ and strictly larger than $3$. It is not $4$ or $6$, so it is $A_4$.

In $S_4$ the nonidentity element orders are $2$, $3$, and $4$. A transposition generates a maximal cyclic subgroup of order $2$: a $4$-cycle squares to a double transposition, not to a transposition. A $3$-cycle generates a maximal cyclic subgroup of order $3$. There are six $4$-cycles and three double transpositions, and the square of a $4$-cycle is a double transposition, so exactly two $4$-cycles square to a given double transposition. Those two elements are inverses and generate the unique cyclic subgroup of order $4$ containing that double transposition. Lemma~\ref{lem:uniquecyclic} therefore gives local counts $2$, $3$, and $4$. Every connection set has average weight at most $4$, and $\rho_{\mathrm{gen}}(S_4)\le 4/24=1/6$.

Equality requires every connection element to have local count $4$. The three double transpositions generate only the normal Klein four-group. Thus an equality witness must contain a $4$-cycle. Let $a=(1\,2\,3\,4)$ and $b=(1\,2\,4\,3)$, and set $S=\{a,a^{-1},b,b^{-1}\}$. These are four elements of order $4$, from two distinct cyclic subgroups. With products read from left to right,
\[
\begin{aligned}
 a^{-1}&=(1\,4\,3\,2),&
 b^{-1}&=(1\,3\,4\,2),\\
 ab^{-1}&=(2\,4\,3),&
 (ab^{-1})a&=(1\,2).
\end{aligned}
\]
The transposition $(1\,2)$ and the $4$-cycle $a$ generate $S_4$: the powers of $a$ already supply $1,2,3,4$ in cyclic order, and conjugation by powers of $a$ produces every adjacent transposition. Hence $S$ generates $S_4$ and attains the bound.

In $A_5$ every nonidentity element has order $2$, $3$, or $5$, so each generates a maximal cyclic subgroup of that order. Lemma~\ref{lem:uniquecyclic} gives local counts $2$, $3$, and $5$. Therefore $\rho_{\mathrm{gen}}(A_5)\le 5/60=1/12$, with equality only for an inverse-closed generating set of elements of order $5$. Let $P$ and $Q$ be distinct cyclic subgroups of order $5$, choose generators $x\in P$ and $y\in Q$, and set $S=\{x,x^{-1},y,y^{-1}\}$. The subgroup $H=\gen{S}$ contains two distinct subgroups of order $5$. There is no element of order $10$, $15$, or $20$ in $A_5$, so $A_5$ has no cyclic subgroup of order $15$ and no subgroup isomorphic to $C_{10}$. A subgroup of order $10$ is consequently dihedral and contains a unique Sylow $5$-subgroup. A subgroup of order $12$ is isomorphic to $A_4$ and has no element of order $5$. A subgroup of order $20$ would have index $3$, and the action on its cosets would give a nontrivial homomorphism $A_5\to S_3$. Its kernel would be a proper nontrivial normal subgroup, contradicting the simplicity of $A_5$. A subgroup of order $30$ would likewise have index $2$ and be normal. Thus $H$ cannot have order $5$, $10$, $12$, $15$, $20$, or $30$, so $H=A_5$. The set $S$ attains the bound.
\end{proof}

\section{Exact computation and verification}\label{sec:computation}
The family formulas are proved exactly, and they were also compared with exhaustive search on small groups. This section describes that computation: how local counts and inversion orbits are enumerated, which groups were examined, and which checks are exhaustive rather than evaluations of the level-set candidates alone.

\subsection{Local counts and inverse-orbit search}
The GAP implementation in Appendix~\ref{app:gap-code} enumerates the elements of $G$ and the cyclic subgroup generated by each element. It evaluates the membership tests defining $c_G(s)$ and partitions the nonidentity elements into inverse orbits. The connected density is the maximum of $\sum_{s\in S}c_G(s)/(|G||S|)$ over generating unions of these orbits.

That maximum is computed by dynamic programming, not by listing level-set candidates. Processing one orbit at a time, the state is the subgroup generated so far together with the cardinality of the chosen set; only the largest weight sum is kept for each pair. Every admissible connection set is a union of orbits, so the optimum is exact for every group the program constructs. The number of reachable subgroups can still grow exponentially with the number of orbits. For the families of Theorems~\ref{thm:dihedralconnected} and~\ref{thm:dicconnected}, the level-set formulas remain the closed evaluations; the program is an independent check of those values.

\subsection{Completed checks}
The GAP source for reproducing these checks is given in Appendix~\ref{app:gap-code}. The checks covered every group of order at most $100$ in the following families: all cyclic groups; all noncyclic abelian groups, written as direct products of cyclic groups of prime-power order; all dihedral groups $D_{2n}$ with $3\leq n\leq 50$; all dicyclic groups $\operatorname{Dic}_n$ with $2\leq n\leq 25$, including the generalized quaternion groups of order at most $64$; and $A_4$, $S_4$, and $A_5$. This is $259$ constructed instances. It is not a census of all isomorphism types of order at most $100$, and no GAP Small Groups computation is claimed.

In every cyclic instance the connected density is $1$. In every elementary abelian instance it equals the proved value $p^{1-r}$. For every dihedral instance with $3\le n\le 50$, the exhaustive optimum agrees with the level-set maximum \eqref{eq:dihedralgeneral}, and with \eqref{eq:dihedralprimepower} precisely when $n$ is a prime power. Every generalized quaternion instance agrees with its closed formula, including $Q_8$, $Q_{16}$, $Q_{32}$, and $Q_{64}$. For $A_4$, $S_4$, and $A_5$ it agrees with Theorem~\ref{thm:alternating}, including the optimal connection-set size $4$. The numerical table below is a sample of these outputs, not a classification.

\begin{table}[htbp]
\centering
\caption{Selected exact connected-density outputs. An optimal size is shown, not necessarily the unique optimal valency.}\label{tab:connected}
\begin{tabular}{lccc}
\toprule
Group & $\rho(G)$ & $\rho_{\mathrm{gen}}(G)$ & Optimal size\\
\midrule
$C_6$ & $1$ & $1$ & $3$\\
$C_2^2$ & $1/2$ & $1/2$ & $2$\\
$C_4\times C_2$ & $3/4$ & $11/20$ & $5$\\
$C_6\times C_2$ & $1$ & $3/4$ & $4$\\
$S_3\cong D_6$ & $1/2$ & $4/9$ & $3$\\
$D_8$ & $1/2$ & $7/16$ & $4$\\
$D_{12}$ & $1/2$ & $5/12$ & $4$\\
$D_{24}$ & $1/2$ & $31/72$ & $6$\\
$Q_8$ & $1$ & $3/5$ & $5$\\
$Q_{16}$ & $1$ & $1/2$ & $5$\\
$Q_{32}$ & $1$ & $33/68$ & $17$\\
$\operatorname{Dic}_3$ & $1$ & $8/15$ & $5$\\
$A_4$ & $1/4$ & $1/4$ & $4$\\
$S_4$ & $1/6$ & $1/6$ & $4$\\
$A_5$ & $1/12$ & $1/12$ & $4$\\
\bottomrule
\end{tabular}
\end{table}

\paragraph{Open problem.}
Let $p$ and $q$ be distinct primes, and let $1\leq u\leq v$ and $a\geq 1$.
Find a closed formula for the connected intersection density
$\rho_{\mathrm{gen}}((C_{p^u}\times C_{p^v})\times C_{q^a})$
in terms of $p,q,u,v,a$.

\appendix
\section{GAP implementation}\label{app:gap-code}
The following program implements the exact inverse-orbit search described
in Section~\ref{sec:computation}. For a nontrivial finite group, it returns
the unrestricted and connected intersection densities, together with all
connection-set sizes that attain the connected maximum.

\begin{lstlisting}[basicstyle=\small\ttfamily,breaklines=true,columns=fullflexible,keepspaces=true,showstringspaces=false,caption={GAP implementation of the exact connected-density algorithm.},label={lst:connected-density-gap}]
# Exact connected intersection density for a nontrivial finite GAP group.
# Read("connected_density.g");; CDExact(AlternatingGroup(4));
CDExact := function(G)
    local el, n, e, cyc, elems, weights, seen, orbits, i, j, s, h, y,
          groups, keys, dp, next, orb, a, w, id, dest, H, k, val,
          full, best, sizes, rho;
    if not IsGroup(G) or not IsFinite(G) or Size(G) < 2 then
        Error("expected a nontrivial finite group");
    fi;
    el := AsSSortedList(G);
    n := Length(el);
    e := One(G);
    cyc := List(el, x -> AsSSortedList(Group([x])));
    weights := List(el, x -> 0);
    seen := List(el, x -> false);
    orbits := [];
    for i in [1..n] do
        s := el[i];
        if s <> e and not seen[i] then
            j := Position(el, s^-1);
            seen[i] := true;
            seen[j] := true;
            w := 0;
            for h in [1..n] do
                y := s * el[h];
                if y in cyc[h] or el[h] in cyc[Position(el, y)] then
                    w := w + 1;
                fi;
            od;
            weights[i] := w;
            weights[j] := w;
            Add(orbits, Set([i,j]));
        fi;
    od;
    rho := Maximum(weights) / n;

    # Each row represents one subgroup. Entry k+1 stores the best
    # sum of local counts for a set of size k generating that subgroup.
    groups := [TrivialSubgroup(G)];
    keys := [[e]];
    dp := [[0]];
    for orb in orbits do
        elems := el{orb};
        a := Length(orb);
        w := Sum(orb, x -> weights[x]);
        next := List(dp, ShallowCopy);  # transitions that omit this orbit
        for id in [1..Length(dp)] do
            H := Group(Concatenation(GeneratorsOfGroup(groups[id]), elems));
            dest := Position(keys, AsSSortedList(H));
            if dest = fail then
                Add(groups, H);
                Add(keys, AsSSortedList(H));
                Add(next, []);
                dest := Length(groups);
            fi;
            for k in [1..Length(dp[id])] do
                if IsBound(dp[id][k]) then
                    val := dp[id][k] + w;
                    if not IsBound(next[dest][k+a]) or
                       val > next[dest][k+a] then
                        next[dest][k+a] := val;
                    fi;
                fi;
            od;
        od;
        dp := next;
    od;
    full := Position(keys, el);
    best := 0;
    sizes := [];
    for k in [1..n-1] do
        if IsBound(dp[full][k+1]) then
            val := dp[full][k+1] / (n*k);
            if val > best then
                best := val;
                sizes := [k];
            elif val = best then
                Add(sizes, k);
            fi;
        fi;
    od;
    return rec(rho := rho, rhoGen := best, optimalSizes := sizes);
end;
\end{lstlisting}


\begin{thebibliography}{99}
\bibitem{aalipour2017}
G.~Aalipour, S.~Akbari, P.~J.~Cameron, R.~Nikandish, and F.~Shaveisi,
On the structure of the power graph and the enhanced power graph of a group,
\emph{Electronic Journal of Combinatorics} \textbf{24} (2017), no.~3, \#P3.16.
\href{https://doi.org/10.37236/6497}{doi:10.37236/6497}.

\bibitem{survey}
J.~Abawajy, A.~Kelarev, and M.~Chowdhury,
Power graphs: a survey,
\emph{Electronic Journal of Graph Theory and Applications} \textbf{1} (2013), no.~2, 125--147.
\href{https://doi.org/10.5614/ejgta.2013.1.2.6}{doi:10.5614/ejgta.2013.1.2.6}.

\bibitem{hosoya}
F.~Ali, B.~A.~Rather, A.~Din, T.~Saeed, and A.~Ullah,
Power graphs of finite groups determined by Hosoya properties,
\emph{Entropy} \textbf{24} (2022), no.~2, article 213.
\href{https://doi.org/10.3390/e24020213}{doi:10.3390/e24020213}.

\bibitem{bera2025}
S.~Bera and H.~K.~Dey,
An exact enumeration of vertex connectivity of the enhanced power graphs of finite nilpotent groups,
\emph{Acta Mathematica Hungarica} \textbf{175} (2025), 550--561.
\href{https://doi.org/10.1007/s10474-025-01524-4}{doi:10.1007/s10474-025-01524-4}.

\bibitem{bubboloni2017quotient}
D.~Bubboloni, M.~A.~Iranmanesh, and S.~M.~Shaker,
Quotient graphs for power graphs,
\emph{Rendiconti del Seminario Matematico della Università di Padova} \textbf{138} (2017), 61--89.
\href{https://doi.org/10.4171/RSMUP/138-3}{doi:10.4171/RSMUP/138-3}.

\bibitem{bubboloni2017alternating}
D.~Bubboloni, M.~A.~Iranmanesh, and S.~M.~Shaker,
On some graphs associated with the finite alternating groups,
\emph{Communications in Algebra} \textbf{45} (2017), no.~12, 5355--5373.
\href{https://doi.org/10.1080/00927872.2017.1307381}{doi:10.1080/00927872.2017.1307381}.

\bibitem{bubboloni2025}
D.~Bubboloni and N.~Pinzauti,
Critical classes of power graphs and reconstruction of directed power graphs,
\emph{Journal of Group Theory} \textbf{28} (2025), no.~3, 713--739.
\href{https://doi.org/10.1515/jgth-2023-0181}{doi:10.1515/jgth-2023-0181}.

\bibitem{cameron}
P.~J.~Cameron,
The power graph of a finite group, II,
\emph{Journal of Group Theory} \textbf{13} (2010), no.~6, 779--783.
\href{https://doi.org/10.1515/JGT.2010.023}{doi:10.1515/JGT.2010.023}.

\bibitem{cameronghosh}
P.~J.~Cameron and S.~Ghosh,
The power graph of a finite group,
\emph{Discrete Mathematics} \textbf{311} (2011), no.~13, 1220--1222.
\href{https://doi.org/10.1016/j.disc.2010.02.011}{doi:10.1016/j.disc.2010.02.011}.

\bibitem{chattopadhyay}
S.~Chattopadhyay and P.~Panigrahi,
Some relations between power graphs and Cayley graphs,
\emph{Journal of the Egyptian Mathematical Society} \textbf{23} (2015), no.~3, 457--462.
\href{https://doi.org/10.1016/j.joems.2015.01.004}{doi:10.1016/j.joems.2015.01.004}.

\bibitem{das2024}
B.~Das, J.~Ghosh, and A.~Kumar,
The isomorphism problem of power graphs and a question of Cameron,
in \emph{44th IARCS Annual Conference on Foundations of Software Technology and Theoretical Computer Science (FSTTCS 2024)},
Leibniz International Proceedings in Informatics, vol.~323 (2024), 20:1--20:23.
\href{https://doi.org/10.4230/LIPIcs.FSTTCS.2024.20}{doi:10.4230/LIPIcs.FSTTCS.2024.20}.

\bibitem{recentsurvey}
A.~Kumar, L.~Selvaganesh, P.~J.~Cameron, and T.~Tamizh Chelvam,
Recent developments on the power graph of finite groups---a survey,
\emph{AKCE International Journal of Graphs and Combinatorics} \textbf{18} (2021), no.~2, 65--94.
\href{https://doi.org/10.1080/09728600.2021.1953359}{doi:10.1080/09728600.2021.1953359}.

\bibitem{lewis2024}
M.~L.~Lewis and C.~Monetta,
Neighbourhoods, connectivity, and diameter of the enhanced power graph,
preprint (2024), \texttt{arXiv:2408.16545}.

\bibitem{ma2022}
X.~Ma, A.~Kelarev, Y.~Lin, and K.~Wang,
A survey on enhanced power graphs of finite groups,
\emph{Electronic Journal of Graph Theory and Applications} \textbf{10} (2022), no.~1.
\href{https://doi.org/10.5614/ejgta.2022.10.1.6}{doi:10.5614/ejgta.2022.10.1.6}.

\bibitem{mirzargar2026}
M.~Mirzargar, S.~Sorgun, and M.~J.~Nadjafi-Arani,
Finite nilpotent groups determined by the enhanced power graph,
\emph{Bulletin of the Malaysian Mathematical Sciences Society} \textbf{49} (2026), no.~6, article 208.
\href{https://doi.org/10.1007/s40840-026-02209-z}{doi:10.1007/s40840-026-02209-z}.

\bibitem{pourgholi2015}
G.~R.~Pourgholi, H.~Yousefi-Azari, and A.~R.~Ashrafi,
The undirected power graph of a finite group,
\emph{Bulletin of the Malaysian Mathematical Sciences Society} \textbf{38} (2015), no.~4, 1517--1525.
\href{https://doi.org/10.1007/s40840-015-0114-4}{doi:10.1007/s40840-015-0114-4}.

\end{thebibliography}
\end{document}